\documentclass[preprint,3p,12pt,pdf]{elsarticle}

\usepackage{lineno,hyperref}
\modulolinenumbers[5]

\usepackage{hyperref}

\usepackage{epstopdf}

\usepackage{tikz}
\usetikzlibrary{arrows}

\usepackage{pst-node}
\usepackage{tikz-cd}

\usepackage[shellescape]{gmp}

\usepackage{mathrsfs}
\usepackage{amssymb}
\usepackage{amsfonts}
\usepackage{latexsym}
\usepackage{amsmath}
\usepackage{xcolor}
\usepackage{wrapfig}
\usepackage{floatflt}

\usepackage{mathtools}
\usepackage{extarrows}

\usepackage{graphicx}

\usepackage{subcaption}

\def\lb{\label}

\newcommand{\er}[1]{\textrm{(\ref{#1})}}

\def\D{\Delta}

\def\Z{{\mathbb Z}}           
    \def\N{{\mathbb N}}   

\let\ge\geqslant                 \let\le\leqslant

                 \def\ts{\times}

\def\el2{\ell^{\,2}}             \def\1{1\!\!1}

\let\ge\geqslant
\let\le\leqslant

\newcommand{\ca}{\begin{cases}}
\newcommand{\ac}{\end{cases}}
\newcommand{\ma}{\begin{pmatrix}}
\newcommand{\am}{\end{pmatrix}}
\def\eq{\begin{equation}}
\def\qe{\end{equation}}
\def\[{\begin{equation}}
\def\]{\end{equation}}

\begin{document}

\begin{frontmatter}

\title{Dynamics of group formation}

\date{\today}

\author
{Anton A. Kutsenko}

\address{Jacobs University, 28759 Bremen, Germany; email: akucenko@gmail.com}

\begin{abstract}
We consider a specific dynamical system of groups formation. It is based simultaneously on a gradient competition between groups and a strong accumulation inside groups. Such a dynamical system demonstrates interesting behavior of densities of emerging groups of $1$, $2$, $3$, ... elements in a steady-state depending on the densities of one-elements groups in a randomly chosen initial state.
\end{abstract}



\end{frontmatter}


{\section{Introduction}\lb{sec1}}

We consider a specific model of group formation without (or with a low) fission. The main dynamics is the competition or repulsion between neighboring groups. The repulsion depends on the sizes or other characteristics of the neighboring groups. After the repulsion, if some groups meet in one common cell then they join together forming one new group. We call this the accumulation or attraction. The process is repeated until it reaches a steady-state. One of the main questions we would like to discuss is the following: for the known initial random distribution of characteristics of the groups, what is the distribution of the characteristics in the steady-state.

Let us provide the general mathematical formulation of the competition-accumulation group formation.
Let $A$ and $B$ be two abelian groups. Denote $C$ and $D$ the abelian groups of functions from $A$ to $B$ and from $A$ to $A$ respectively. Note that the introduced abelian groups and emerging groups in the model are different notions. Let $T(0)=\{T(0,m)\}_{m\in A}\in C$ be some initial state. Consider the evolution equation
\[\lb{001}
 T(n+1,m)=\sum_{k\in A:\ k-m=(\nabla T(n))(k)}T(n,k),\ \ \ \ \ n\in\N\cup\{0\},\ m\in A,
\]
where $T(n)=\{T(n,m)\}_{m\in A}\in C$, and
\[\lb{002}
 \nabla: C\to D
\]
is some mapping. The group $A$ is the space of the model. The group $B$ describes some characteristics of the model. If $B=\Z$ then $T(n,m)$ can be interpreted as the number of elements in cell $m$ at time $n$. The mapping $\nabla$ shows how neighbors of $T(n,k)$ affect on its moving from the position $k$ to the new position $k-(\nabla T(n))(k)$. For some models, it is reasonable to assume that $\nabla$ is an analog of the standard gradient operator. All the quantities $T(n,k)$ which move to the same cell $m$ will accumulate forming the new quantity $T(n+1,m)$. In some sense, $\{T(n)\}_{n\ge0}$ is a non-linear Markov chain, since the linear operator that maps $T(n)$ to $T(n+1)$ depends on $T(n)$.  Some information about non-linear Markov processes is available in, e.g., \cite{F}. Note that there is some analogy with material derivatives $u\cdot\nabla u$ appearing in Euler and Navier-Stokes type equations, where the linear operator $u\cdot\nabla*$ depends on the argument $u$. This may indicate that nonlinear Markov chains can be as difficult to analyze as the problems of fluid dynamics in some cases. While the model \er{001} is deterministic, we assume further that the initial input $T(0)$ is a random field. This assumption leads to some interesting observations in stochastic analysis of steady-states. The evolution equation \er{001} may describe some non-standard particle dynamics, or a process of teams formation in, e.g., scientific communities, or an evolution of economic and biological systems under specific conditions. I cannot find models identical to \er{001} in the base of knowledge. Some continuous equations of a group formation are considered in \cite{GL}. However, the discrete nature of, e.g., steady-states solutions of \er{001} cannot be reflected completely in continuous equations. Nevertheless, continuous equations are successfully applied in practice, see \cite{DLP}. Concerning discrete dynamical systems, it is useful to note \cite{HSGN}, where some aspects of the game theory are applied for the analysis of evolutionary dynamics experimental data. Finally, let us note \cite{CK} that motivates the current research in some way. And of course, we remember the great ideas that inspire similar research at all times: the Conway's Game of Life and the Turing Machine.

As already mentioned above, one of the main questions is the statistical analysis of the dynamical system \er{001}. Namely, having a big set of randomly chosen initial states $T(0)$ what can we say about the steady-states $T(n)$. In the next section, we consider one- and two-dimensional models, which already demonstrates some interesting statistical properties.

{\section{Main results}\lb{sec2}}
\begin{figure}[h]
\center{\includegraphics[width=0.99\linewidth]{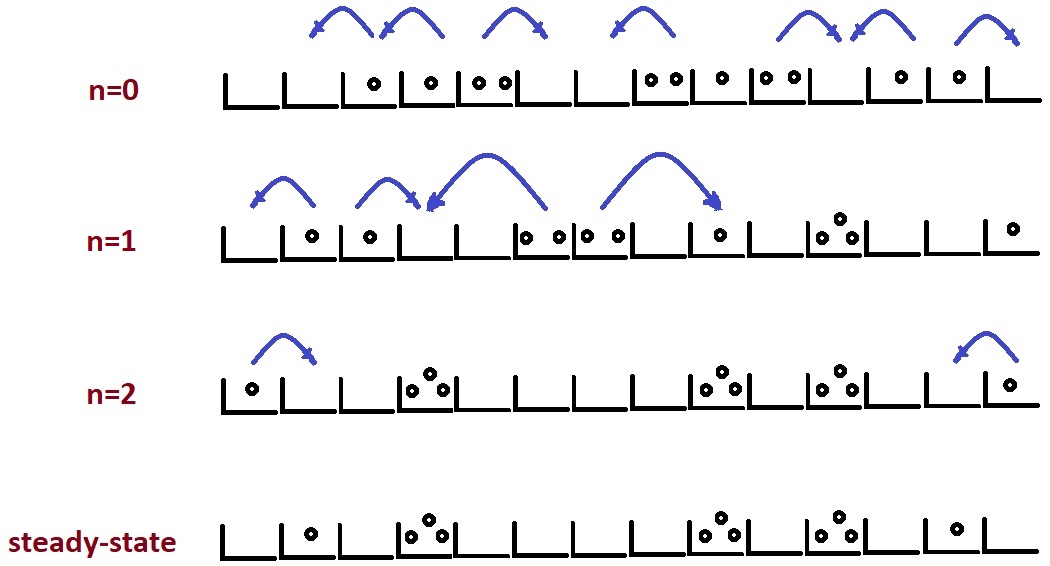}}
\caption{An example of evolution of the system $T(0,0)=T(0,1)=0$, $T(0,2)=T(0,3)=1$, $T(0,4)=2$, $T(0,5)=T(0,6)=0$, $T(0,7)=2$, $T(0,8)=1$, $T(0,9)=2$, $T(0,10)=0$, $T(0,11)=T(0,12)=1$, $T(0,13)=0$ on the torus $\Z_{14}$. }\lb{fig0}
\end{figure}
\subsection{1D model} The space of the model is discrete torus $\Z_M$. At time $n\in\N\cup\{0\}$, each cell $m\in\Z_M$ of the torus contains the group of elements of the size $T(n,m)\in\Z$. If $T(n,m)=0$ then the cell $m$ is empty at time $n$. The evolution of the group formation is given by the equation
\[\lb{101}
 T(n+1,m)=\sum_{k-m=T(n,k+1)-T(n,k-1)}T(n,k),\ \ \ for \ \ \ n\ge0,\ m\in\Z_M.
\]
The equality below the sum symbol assumes by modulo $M$. Each neighbor $T(n,k-1)$ and $T(n,k+1)$ pushes the group $T(n,k)$ away from itself at a distance proportional to a size of the pushing neighbor. After this, the groups occupying the same cell will merge into one group of size equal to the sum of the sizes of the merged groups. The situation is repeated. Very often it reaches the steady-state when each emerging group has no neighbors. One such example is demonstrated in Fig. \ref{fig0}.

For the numerical simulations, we chose $M=3000$ and $M=4000$. The initial state $\{T(0,m)\}_{m\in\Z_M}$ consists of one-element groups and empty cells only. The density of one-element groups is 
$$
 p=\frac{\#\{m:\ T(0,m)=1\}}{M},
$$
where $\#$ denotes the number of elements of the set. In fact, it is assumed that the initial states consist of $M$ independent equally distributed random variables $T(0,m)$, $m\in\Z_M$, which take the value $1$ with the probability $p$ and the value $0$ with the probability $1-p$. For each $p\in[0,0.96]$ with the step $\D p=0.01$ we produce $10000$ random initial states that evolve to a constant steady-state. The steady-state consists of groups of different sizes. Their densities are given by
$$
 Q_r(p)=\frac{\#\{m:\ T({\rm steady-state,}m)=r\}}{M}.
$$ 
For $r=1,2,3,4$ the average densities over $10000$ initial random states are plotted in Fig. \ref{fig1}. The first note is that $Q(p)$ are almost independent on $M$. In fact, we found that $Q_r(p)$ for $M=300$ and $M=3000$ are almost identical if $p$ is not very close to $1$ and $r$ is not large. It is expected, since the evolution equation \er{101} seems homogeneous and the values in any initial state are independent. However, the average time $N_{\rm st}$ of reaching steady-states is different, e.g., for $p=0.8$, we have $N_{\rm st}=50$ for $M=3000$ and $N_{\rm st}=52.5$ for $M=4000$. The plot of the average $N_{\rm st}$ over $10000$ samples for $p=0.35$, $0.6$, and $0.85$ and for the model sizes $300\le M\le 10000$ is presented in Fig. \ref{fig2}. If the difference between model sizes $\D M$ is not large then the difference $\D N_{\rm st}$ is small but it is enough for the formation of large-size groups. Obviously, the large-size groups may appear in large models only. Nevertheless, the density of small-size groups is sufficiently stable starting from small $M$. For small $p$, the initial state is already close to a steady-state since the distribution  of one-element groups is sparse enough and most of them have no neighbors. For $p>0.5$, the density of one-element groups decreases and the two-element groups become dominant for $p\ge0.7$. There is some competition between two- and three-element groups for $p\in[0.7,0.85]$, but two-element groups always lead. As already mentioned above, if $p$ is close to $1$ then large groups will appear with, perhaps, interesting dynamics involving competitions between large groups and maybe a change of leaders among them. 

We found also that the probability of periodic states (steady-states, non constant in time) for $p\approx1$ increases. At the same time, for $p<0.96$ and $M>300$ almost all initial states reach constant steady-states.

\begin{figure}[h]
    \centering
    \begin{subfigure}[b]{0.49\textwidth}
        \includegraphics[width=\textwidth]{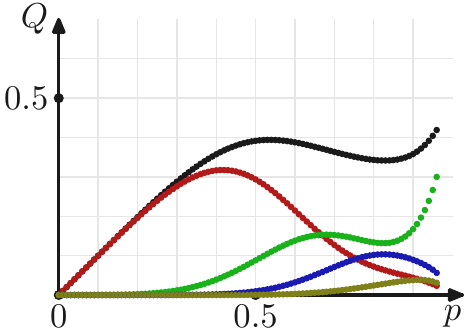}
        \caption{M=3000}
    \end{subfigure}
    \begin{subfigure}[b]{0.49\textwidth}
        \includegraphics[width=\textwidth]{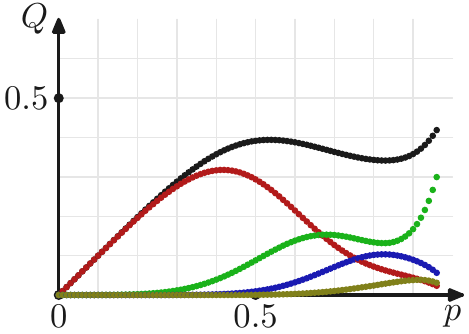}
        \caption{M=4000}
    \end{subfigure}
    \caption{For the initial density $p$ of one-element groups, the steady-state densities $Q_r(p)$ of emerging groups of $r=1$ (red), $r=2$ (green), $r=3$ (blue), and $r=4$ (dark green) elements in 1D model are presented. Black dots show the total density of all groups. The parameter $M$ is the size of the model.}\lb{fig1}
\end{figure}

\begin{figure}[h]
	\center{\includegraphics[width=0.5\linewidth]{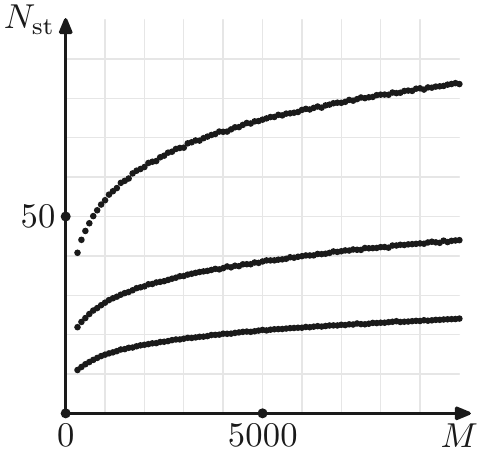}}
	\caption{The average time $N_{\rm st}$ of reaching steady-state is computed for different model sizes $M$. The bottom curve corresponds to $p=0.35$, the middle curve corresponds to $p=0.6$, and the top curve corresponds to $p=0.85$, where $p$ is the probability of one-element groups in the initial state. }\lb{fig2}
\end{figure}

\subsection{Primitive one-step 1D model}
While considered 1D model looks simple, it is sufficiently complex to derive analytic formulas for the densities of emerging groups $Q_r(p)$ in steady-states, and for the relaxation time $N_{\rm st}$. Let us consider very simplified model for which we can derive some analytic formulas. We consider $\Z_M$, where independent equally distributed random variables $\zeta_m$ are placed at each cell $m\in\Z_M$. The variables $\zeta_m$ take two values $0$ and $1$ with the probabilities $1-p$ and $p$ respectively, with $p\in[0,1]$. The variables represent the random distribution of on-element groups. Then, say each even group moves left or right with the probability $\frac12$ and merge with the group in the corresponding neighboring cell. After this one step, the system becomes in steady-state since all the neighboring cells of each odd cell is empty. Here, we assume that $M$ is even. It is easy to see that in the steady-sate each odd cell may contain one, two, or three random variables with the probabilities $s_1=\frac14$, $s_2=\frac12$, and $s_3=\frac14$ respectively. Using this and the fact that
$$
 {\rm Prob}(\zeta_a=0)=1-p,\ \ \ {\rm Prob}(\zeta_a=1)=p,
$$
$$
 {\rm Prob}(\zeta_a+\zeta_b=0)=(1-p)^2,\ \ \ {\rm Prob}(\zeta_a+\zeta_b=1)=2p(1-p),\ \ \ {\rm Prob}(\zeta_a+\zeta_b=2)=p^2,
$$
\begin{multline}
 {\rm Prob}(\zeta_a+\zeta_b+\zeta_c=0)=(1-p)^3,\ \ \ {\rm Prob}(\zeta_a+\zeta_b+\zeta_c=1)=3p(1-p)^2,\\
 {\rm Prob}(\zeta_a+\zeta_b+\zeta_c=2)=3p^2(1-p),\ \ \ {\rm Prob}(\zeta_a+\zeta_b+\zeta_c=3)=p^3\notag
\end{multline}
for independent equally distributed two-valued random variables  $\zeta_a$, $\zeta_b$, and $\zeta_c$ introduced above, we deduce that the densities $Q_r(p)$ of $r$-element groups for $r=1,2,3$ are
\begin{multline}\lb{200}
 Q_r(p)=\frac{s_1{\rm Prob}(\zeta_a=r)\frac{M}2+s_2{\rm Prob}(\zeta_a+\zeta_b=r)\frac{M}2+s_3{\rm Prob}(\zeta_a+\zeta_b+\zeta_c=r)\frac{M}2}M=\\
 \ca
  \frac{8p-10p^2+3p^3}8,& r=1,\\
  \frac{5p^2-3p^3}8,& r=2,\\
  \frac{p^3}8,& r=3.
 \ac
\end{multline}
\begin{figure}[h]
	\center{\includegraphics[width=0.5\linewidth]{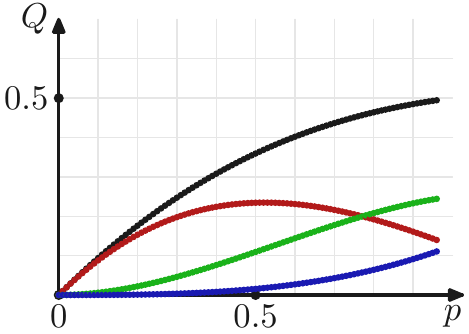}}
	\caption{For the initial density $p$ of one-element groups, the steady-state densities $Q_r(p)$, see \er{200}, of emerging groups of $r=1$ (red), $r=2$ (green), and $r=3$ (blue) elements in the primitive model are presented. Black dots show the total density of all groups. }\lb{fig3}
\end{figure}
The densities $Q_r(p)$ for the primitive model are plotted in Fig. \ref{fig3}. While the curves are much simpler than the similar ones from Fig. \ref{fig1}, they demonstrate the same general tendencies: the density of one-element groups decreases and the two-element groups are leaders in dense systems. 

\subsection{2D model}
\begin{figure}[h]
    \centering
    \begin{subfigure}[b]{0.49\textwidth}
        \includegraphics[width=\textwidth]{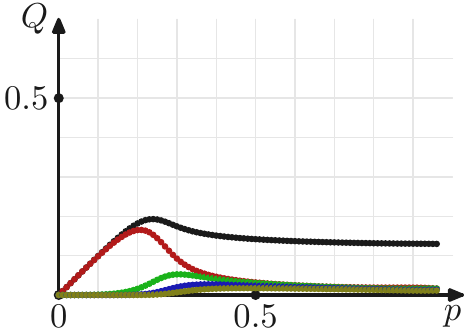}
        \caption{$M_1=100,\ M_2=150$}
    \end{subfigure}
    \begin{subfigure}[b]{0.49\textwidth}
        \includegraphics[width=\textwidth]{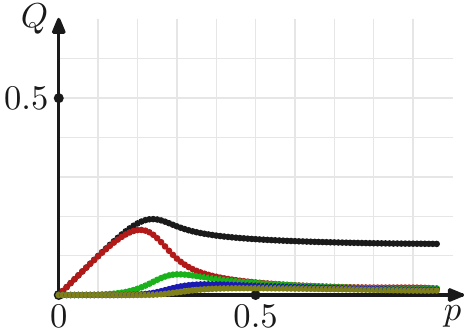}
        \caption{$M_1=M_2=200$}
    \end{subfigure}
    \caption{For the initial density $p$ of one-element groups, the steady-state densities $Q_r(p)$ of emerging groups of $r=1$ (red), $r=2$ (green), $r=3$ (blue), and $r=4$ (dark green) elements in 2D model are presented. Black dots show the total density of all groups. The parameters $M_1$, $M_2$ show the size of the model.}\lb{fig4}
\end{figure}
\begin{figure}[h]
	\centering
    \begin{subfigure}[b]{0.49\textwidth}
        \includegraphics[width=\textwidth]{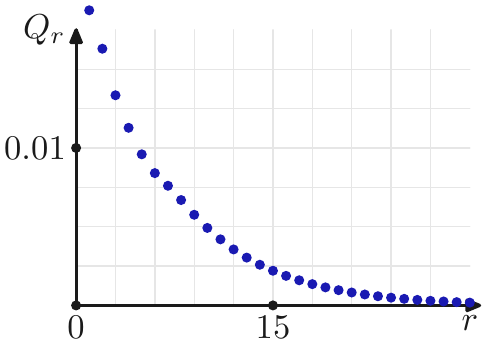}
        \caption{$p=0.9$}
    \end{subfigure}
    \begin{subfigure}[b]{0.49\textwidth}
        \includegraphics[width=\textwidth]{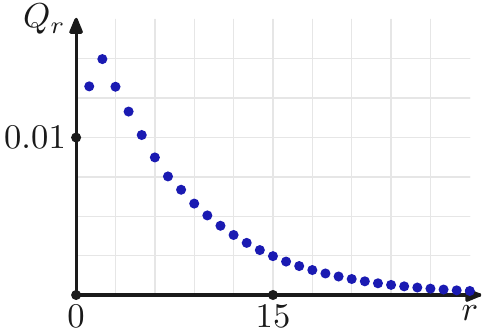}
        \caption{$p=0.99$}
    \end{subfigure}        
        \caption{The steady-state densities $Q_r(p)$ of $r$-element groups in 2D model of size $M_1=M_2=200$ are presented. }\lb{fig5}
\end{figure}
Finally, let us briefly consider 2D model $\Z_{M_1}\ts\Z_{M_2}$ with the evolution equation
\[\lb{300}
 T(n+1,{\bf m})=\sum_{{\bf k}-{\bf m}=\nabla T(n,{\bf k})}T(n,{\bf k}),\ \ for\ \ n\ge0,\ \ {\bf m}\in\Z_{M_1}\ts\Z_{M_2},
\]
where
\[\lb{301}
 \nabla T(n,{\bf k})=\ma T(n,k_1+1,k_2)-T(n,k_1-1,k_2) \\
  T(n,k_1,k_2+1)-T(n,k_1,k_2-1) \am,\ \ \ {\bf k}=\ma k_1 \\ k_2 \am\in\Z_{M_1}\ts\Z_{M_2}.
\]
As in 1D model, we generate a lot of random one-element groups with densities $p\in[0,0.96]$ and compute the average densities $Q_r(p)$ of emerging $r$-element groups in steady-states. These curves are plotted in Fig. \ref{fig4}. While there are some similarities with Fig. \ref{fig1}, the curves are different. Again, the dependence $Q_r(p)$ on $M_1$ and $M_2$ is small if $M_1$ and $M_2$ are sufficiently large. Also, the difference between the average densities of $1$-, $2$-, $3$-, $4$-element groups in 2D model is not large in comparison with 1D model. The average densities $Q_r(p)$ for $p=0.9$, $p=0.99$ and for the model size $M_1=M_2=200$ are presented in Fig. \ref{fig5}.

{\section{Conclusion}\lb{sec3}} The statistical analysis of the dynamical model of group formation with an out-group repulsion and a strong in-group attraction is provided. Depending on the probability of one-element groups in the initial state, the distribution of groups of different sizes in steady-states is sufficiently complex. However, two-element groups demonstrate domination in dense systems. We also found a primitive analytic model that can emulate some of the densities tendencies. At the same, the derivation of analytic results for the densities and for the relaxation time in the considered non-primitive models seems quite hard.

\clearpage

\section*{Acknowledgements} 
This paper is a contribution to the project M3 of the Collaborative Research Centre TRR 181 "Energy Transfer in Atmosphere and Ocean" funded by the Deutsche Forschungsgemeinschaft (DFG, German Research Foundation) - Projektnummer 274762653. This work is also supported by the RFBR (RFFI) grant No. 19-01-00094.


\begin{thebibliography}{9}
\bibitem{F}
T. D. Frank,  "Markov chains of nonlinear Markov processes and an application to a winner-takes-all model for social conformity". {\it J. Phys. A} {\bf 41}, 282001, 2008.

\bibitem{GL}
S. Gueron and S. A. Levin,  "The dynamics of group formation". {\it Math. Biosci.} {\bf 128}, 243-264, 1995.

\bibitem{DLP}
P. Degond, J.-G. Liu, and R. L. Pego,  "Coagulation–fragmentation model for animal group-size statistics". {\it J. Nonlinear. Sci.} {\bf 27}, 379-424, 2017.

\bibitem{HSGN}
M. Hoffman, S. Suetens, U. Gneezy, and M. A. Nowak,  "An experimental investigation of evolutionary dynamics in the Rock-Paper-Scissors game". {\it Sci. Rep.} {\bf 5}, 8817, 2015.

\bibitem{CK}
F. F. Chen and D. T. Kenrick,  "Repulsion or attraction? Group membership and assumed attitude similarity". {\it J. Pers. Soc. Psychol.} {\bf 83}, 111-125, 2002.

\end{thebibliography}
\end{document}